\documentclass[12pt,a4paper]{amsart}

\theoremstyle{plain}
\newtheorem{lem}{Lemma}
\newtheorem*{thm}{Theorem}
\newtheorem*{cor}{Corollary}
\theoremstyle{definition}
\newtheorem{rmk}{Remark}

\newcommand{\itr}[2]{#1^{\circ #2}}
\newcommand{\CC}{\mathbb{C}}
\newcommand{\NN}{\mathbb{N}}
\newcommand{\ie}{i.\,e. }
\newcommand{\mthm}{the main theorem}

\title{Minimal polynomial of an exponential automorphism of $\mathbb{C}^n$}
\author[J. Zygad{\l}o]{Jakub Zygad{\l}o}
\date{January 4, 2008}
\address{Institute of Mathematics\\
Jagiellonian University\\
Reymonta 4\\
30-059 Krak\'ow, Poland.}
\email{jakub.zygadlo@im.uj.edu.pl}

\begin{document}
\begin{abstract}
We show that the minimal polynomial of a polynomial exponential automorphism $F$ of $\mathbb{C}^n$ (i.\,e. $F=\exp(D)$ where $D$ is a locally nilpotent derivation) is of the form $\mu_F(T)=(T-1)^d$, $d=\min\{m\in\mathbb{N}: D^{\circ m}(X_i)=0\text{ for } i=1,\ldots,n\}$. 
\end{abstract}
\maketitle

\section{Introduction}

Let $k$ be a field of characteristic zero and let $A$ be a $k$-algebra. Recall that a {\it $k$-derivation} of $A$ is a $k$-linear mapping $D\colon A\to A$ fulfilling the Leibniz rule $D(ab)=D(a)b+aD(b)$. We will write $\itr{D}{n}$ for the $n$-th iterate of $D$, \ie $\itr{D}{n}=D\circ\itr{D}{(n-1)}$ and $\itr{D}{0}=I$ - the identity. If for every $a\in A$ there exists $n=n(a)\in\NN$ such that $\itr{D}{n}(a)=0$, derivation $D$ is called {\it locally nilpotent}.

If $D$ is a locally nilpotent derivation of $A$, we define the {\it exponential} of $D$, denoted $\exp(D)$, by the formula
$$\exp(D)(a):=\sum_{i=0}^{\infty}\frac{1}{i!}\itr{D}{i}(a)$$
It is easy to see that $\exp(D)\colon A\to A$ is a $k$-endomorphism of $A$. One can also check that if locally nilpotent $k$-derivations $D$ and $E$ commute (\ie $D\circ E=E\circ D$), then $\exp(D)\circ\exp(E)= \exp(E)\circ\exp(D)= \exp(D+E)$. Therefore, $\exp(D)$ is an automorphism of $A$ with the inverse $\exp(D)^{-1}=\exp(-D)$.
In the paper we prove the following

\begin{thm}
Let $D$ be a locally nilpotent derivation of $\CC[X_1,\ldots,X_n]$, $F:=(\exp(D)(X_1),\ldots,\exp(D)(X_n)\,)\colon\CC^n \to \CC^n$ and  $d:=\min\{m\in\nolinebreak\NN: \itr{D}{m}(X_i)=0\textnormal{ for }i=1,\ldots,n\}$. Then the minimal polynomial for $F$ equals $\mu_F(T)=(T-1)^d=\sum_{j=0}^d(-1)^{d-j}\binom{d}{j}T^j$ (\ie the mapping $\mu_F(F)=\sum_{j=0}^d(-1)^{d-j}\binom{d}{j}\itr{F}{j}$ is zero and $p(F)\neq 0$ for any polynomial $p\in \CC[T]\setminus\nolinebreak\{0\}$ of degree less than $d$).
In particular, we have the following formula for the inverse of $F$:
$$F^{-1}=\sum_{j=0}^{d-1}(-1)^j\binom{d}{j+1}\itr{F}{j}$$
\end{thm}

\section{Preparatory steps}
\noindent Firstly, we will prove two simple lemmas:

\begin{lem}\label{derlem}
Let $A$ be a $k$-algebra, $D$ - a locally nilpotent $k$-derivation of $A$ and $a\in A$. If for some $m\geq 1$ and $\alpha_0,\ldots,\alpha_{m-1}\in k$ there is an equality $$\itr{D}{m}(a)=\sum_{i=0}^{m-1}\alpha_i\itr{D}{i}(a),$$
then $\itr{D}{m}(a)=0$.
\begin{proof}
We will proceed by induction on $m$. If $m=1$ we have $D(a)=\alpha_0 a$ and the result is well known (even for $\alpha_0\in A$, if $A$ has no zero divisors - see for example \cite{Es}, Prop. 1.3.32), but we will prove it for the sake of completeness. If $D(a)=\alpha_0 a$, then $\itr{D}{n}(a)=\itr{D}{(n-1)}(\alpha_0 a)=\ldots=D(\alpha_0^{n-1}a)=\alpha_0^n a$ for all $n\in\NN$. Because $D$ is locally nilpotent, we must have $\itr{D}{n}(a)=0$ for some $n$ and consequently $\alpha_0=0$ or $a=0$. Now let $m>1$ and assume that the lemma holds for all $m'<m$. Suppose $\itr{D}{m}(a)\neq 0$ and let $M\in\NN$ be such that $\itr{D}{M}(a)=0$ and $\itr{D}{(M-1)}(a)\neq 0$ (note $M>m$). Set $i_0:=\max\{0\leq i<m:\alpha_i\neq 0\}$, so we can write $0=\itr{D}{M}(a)=\itr{D}{(M-m)}(\itr{D}{m}(a))= \itr{D}{(M-m)}\big(\sum_{i=0}^{i_0}\alpha_i\itr{D}{i}(a)\big)= \sum_{i=0}^{i_0}\alpha_i\itr{D}{i}\big(\itr{D}{(M-m)}(a)\big)$. Let $a':=\itr{D}{(M-m)}(a)$. Because $\alpha_{i_0}\neq 0$, we have $\itr{D}{i_0}(a')=-\sum_{i=0}^{i_0-1}\frac{\alpha_i}{\alpha_{i_0}}\itr{D}{i}(a')$ and since $i_0<m$, we obtain $\itr{D}{i_0}(a')=0$ by the induction hypothesis - this is a contradiction with $\itr{D}{i_0}(a')=\itr{D}{(M-m+i_0)}(a)\neq 0$.
\end{proof}
\end{lem}

\begin{lem}\label{coeflem}
Let $d>0$, $i\in\NN$ and define $$\beta_{d,i}:=\sum_{m=0}^d(-1)^m\binom{d}{m}m^i$$
We have $\beta_{d,i}=0$ if and only if $i<d$.
\begin{proof}
Equality $\beta_{d,0}=0$ follows from expansion of $(1-1)^d=0$ and the case $d=1$ is obvious. Let $d>1$, $i>0$ and proceed by induction on $d$. We have
\small
\begin{align*}
\beta_{d,i}&=\sum_{m=1}^d(-1)^md\binom{d-1}{m-1}m^{i-1}= 
-d\sum_{m=0}^{d-1}(-1)^m\binom{d-1}{m}(m+1)^{i-1}=\\ 
&=-d\sum_{j=0}^{i-1}\binom{i-1}{j}\left(\sum_{m=0}^{d-1}(-1)^m\binom{d-1}{m}m^j\right)= -d\sum_{j=0}^{i-1}\binom{i-1}{j}\beta_{d-1,j}
\end{align*}
\normalsize
and for $i<d$ we conclude by the induction hypothesis, because all $\beta_{d-1,j}=0$. To deal with the case $i\geq d$, note that $\beta_{1,i}=-1$ for $i\geq 1$, so $\beta_{2,i}=-d\sum_{j=0}^{i-1}\binom{i-1}{j}\beta_{1,j}>0$ for $i\geq 2$. Proceeding in this way, we see that $(-1)^d\beta_{d,i}>0$ for $i\geq d$.
\end{proof}
\end{lem}

From now on we will focus our attention on the case $k=\CC$ and $A=\CC[X_1,\ldots,X_n]$ - the ring of polynomials in $n$ variables.
It can be shown that every $\CC$-derivation $D$ of $A$ is of the form $D=\sum_{i=1}^nf_i\partial_{x_i}$ for some $f_1,\ldots,f_n\in A$, where $\partial_{x_i}=\frac{\partial}{\partial X_i}$ is the standard differential with respect to $X_i$.

If $\Phi\colon A\to A$ is a $\CC$-endomorphism of $A$, one can define a polynomial mapping $\Phi_*\colon\CC^n\to\CC^n$ by
$$\Phi_*=(\Phi(X_1),\ldots,\Phi(X_n))$$
Obviously $I_*=I$ and $(\Phi\circ\Psi)_*=\Psi_*\circ\Phi_*$, so each $\CC$-automorphism $\Phi$ of $A$ gives rise to a polynomial automorphism $\Phi_*$ of the affine space $\CC^n$. In particular, if $D$ is a locally nilpotent derivation of $A$ and $\Phi=\exp(D)$, we have an automorphism $F=\exp(D)_*=(\exp(D)(X_1),\ldots,\exp(D)(X_n))$ of $\CC^n$, called the {\it exponential automorphism}.

In \cite{FM}, the following class of polynomial automorphisms is considered:
Let $F=(F_1,\ldots,F_n)$ be a polynomial automorphism of $\CC^n$. If there is an univariate polynomial $p(T)\in\CC[T]\setminus\{0\}$ such that $p(F)=0$ (\ie if $p(T)=a_0+a_1T+\ldots+a_mT^m$ this means $a_0I+a_1F+\ldots+a_m\itr{F}{m}=0$), then $F$ is called {\it locally finite}.

It is easy to see that the set $I_F:=\{p\in\CC[T]: p(F)=0\}$ forms an ideal in $\CC[T]$; its monic generator will be called {\it minimal polynomial} for $F$ and denoted $\mu_F$. The paper \cite{FM} gives many equivalent conditions for $F$ to be locally finite and a formula for a polynomial $p(T)$ such that $p(F)=0$, provided $F(0)=0$ (see \cite{FM}, Th. 1.2). Unfortunately, there is no such result when $F(0)\neq 0$ and it is not easy to find the minimal polynomial $\mu_F$, either. We solve this problem for exponential automorphisms of $\CC^n$ in the following section.

\section{Main result and its consequences}

\begin{thm}[main theorem]
Let $D$ be a locally nilpotent derivation of $\CC[X_1,\ldots,X_n]$, $F:=\exp(D)_*$ and $d:=\min\{m\in\NN: \itr{D}{m}(X_i)=\nolinebreak 0 \text{ for } i=1,\ldots,n\}$. Then the minimal polynomial for $F$ equals $\mu_F(T)=(T-1)^d$. 
\begin{proof}
Note that for $m\in\NN$, we have $\itr{F}{m}=(\itr{\exp(D)}{m})_*=\exp(mD)_*$ (because $D$ commutes with $D$), so if $F=(F_1,\ldots,F_n)$ then $$(\itr{F}{m})_j= \sum_{i=0}^{d-1}\frac{1}{i!}\itr{(mD)}{i}(X_j)=  \sum_{i=0}^{d-1}\frac{1}{i!}m^i\itr{D}{i}(X_j), \quad j=1,\ldots,n$$
Since
\footnotesize
$$\sum_{m=0}^d(-1)^m\binom{d}{m}(\itr{F}{m})_j=  \sum_{i=0}^{d-1}\frac{1}{i!}\left(\sum_{m=0}^d(-1)^m\binom{d}{m}m^i\right)\itr{D}{i}(X_j)=  \sum_{i=0}^{d-1}\frac{1}{i!}\beta_{d,i}\itr{D}{i}(X_j)$$
\normalsize
we conclude by Lemma \ref{coeflem} that $\sum_{m=0}^d(-1)^m\binom{d}{m}\itr{F}{m}=0$.
This argument shows that the polynomial $(1-T)^d\in I_F=\{p\in\CC[T]: p(F)=0\}$. To prove minimality of its degree, assume for example $d=\min\{m\in\NN: \itr{D}{m}(X_1)=0\}$ and suppose that $\mu_F(T)=(T-1)^e$ for some $e<d$. Then $0=(-1)^e\big(\mu_F(F)\big)_1= \sum_{m=0}^{e}(-1)^m\binom{e}{m}(\itr{F}{m})_1= \sum_{i=0}^{d-1}\frac{1}{i!}\beta_{e,i}\itr{D}{i}(X_1)$ and $\beta_{e,d-1}\neq 0$ by Lemma \ref{coeflem}. Therefore $\itr{D}{(d-1)}(X_1)= -\sum_{i=0}^{d-2}\frac{1}{i!}\frac{\beta_{e,i}}{\beta_{e,d-1}}\itr{D}{i}(X_1)$ and due to Lemma \ref{derlem} we get $\itr{D}{(d-1)}(X_1)=0$, despite the definition of $d$ - a contradiction.
\end{proof}
\end{thm}

\begin{cor}
Since $\mu_F(F)=0$, we have $I=\Big(\sum\limits_{m=1}^d(-1)^{m-1}\binom{d}{m}\itr{F}{(m-1)}\Big)\circ\nolinebreak F$ and therefore the inverse of $F$ is given by
\small
$$F^{-1}=\sum_{m=0}^{d-1}(-1)^m\binom{d}{m+1}\itr{F}{m}$$
\normalsize
\end{cor}

\begin{rmk}
The famous Nagata automorphism of $\CC^3$ (see \cite{Na}) defined by $N=(X-2Y\sigma-Z\sigma^2,Y+Z\sigma,Z)$ where $\sigma=XZ+Y^2$ can be seen as an exponential of a locally nilpotent derivation $D=-2Y\sigma\partial_{x}+Z\sigma\partial_{y}$ of $\CC[X,Y,Z]$. It is easy to check that $D(\sigma)=0$ and $\itr{D}{3}(X)=\itr{D}{3}(Y)=\itr{D}{3}(Z)=0$, so \mthm\ gives $\mu_N(T)=(T-1)^3$, whereas by (\cite{FM}, Th. 1.2) we only get that $p(T)=(T-1)^{55}\in I_N$.
\end{rmk}

\begin{rmk}
Let $d\geq 2$ and $D=Y^{d-2}\partial_x+\partial_y$. Obviously $D$ is a locally nilpotent derivation of $\CC[X,Y]$ and $\itr{D}{d}(X)=\itr{D}{d}(Y)=0$ ($d$ is minimal). If we let $F=\exp(D)_*$, then $\mu_F(T)=(T-1)^d$ by \mthm. Since clearly $\deg F=d-2$, this shows that the estimate $\deg\mu_F\leq\deg F+1$ (\cite{FM}, Th. 4.2) need not hold if $F(0)\neq 0$.
\end{rmk}

\begin{rmk}
Recall that if $P=(P_1,\ldots,P_n)\colon\CC^n\to\CC^n$ is a polynomial mapping, then $P^*$ given by $P^*(X_i):=P_i\in\CC[X_1,\ldots,X_n]$ defines a $\CC$-endomorphism of $\CC[X_1,\ldots,X_n]$.
Let $F=(X+g(Y,Z),Y+h(Z),Z)$ be an upper triangular automorphism of $\CC^3$ ($g\in\CC[Y,Z]$, $h\in\CC[Z]$). If $g=0$ or $h=0$, then $F$ is easily seen to be an exponential of a locally nilpotent derivation of $\CC[X,Y,Z]$. So let us suppose that $g\neq 0$ and $h\neq 0$. We will show that the minimal polynomial for $F$ equals $\mu_F(T)=(T-1)^d$, where $d:=2+\deg_Y\!g$ and therefore in this case we also have $F=\exp(D)_*$ (see \cite{FM}, Th. 2.3) for the locally nilpotent derivation $D$ given by the following formula: $D=\sum_{m=1}^{+\infty}\frac{(-1)^{m+1}}{m}\itr{(F^*-\nolinebreak I^*)}{m}$ (cf. \cite{Es}, Ch. 2). Obviously $D(Z)=0$, $D(Y)=h(Z)$ and one can use above formula to evaluate $D(X)$ - note that if we can show that the minimal polynomial has degree $d$, then only first $d-1$ summands are nonzero.
Write $g(Y,Z)=\sum_{i=0}^{d-2} Y^ig_i(Z)$. Iterating $F$, we get
\small
\begin{align*}
\itr{F}{m}&=(X+\sum_{j=0}^{m-1}\sum_{i=0}^{d-2}\big(Y+jh(Z)\big)^ig_i(Z),Y+mh(Z),Z)=\\
&=(X+\sum_{i=0}^{d-2}g_i(Z)\sum_{k=0}^i\binom{i}{k}Y^{i-k}h(Z)^k\sum_{j=0}^{m-1}j^k,Y+mh(Z),Z)
\end{align*}
\normalsize
Let $s_k(m):=\sum_{j=0}^{m-1}j^k$ and note that $s_k$ is a polynomial in $m$ of degree $k+1\leq d-1<d$. Therefore, Lemma \ref{coeflem} gives $\sum_{m=0}^{d}(-1)^m\binom{d}{m}s_k(m)=\nolinebreak 0$ for all $k\leq d-2$ and we can argue as in the proof of \mthm\ (since $g_{d-2}(Z)h(Z)^{d-2}s_{d-2}(m)$ is the only term involving $m^{d-1}$, we must have $\itr{D}{(d-1)}(X)\neq 0$). Consequently $\mu_F(T)=(T-1)^d$ and $F=\exp(D)_*$, where 
\small
$$D=\bigg(\sum_{m=1}^{d-1}\sum_{i=1}^{m}\frac{(-1)^{i+1}}{m}\binom{m}{i}\big((\itr{F}{i})_1-X\big)\bigg)\partial_x+h(Z)\partial_y$$ 
\normalsize
and $d=2+\deg_Y\!g$. Note that $d$ is minimal and easily found in this case (there are obstacles to calculations of the minimal degree, cf. \cite{FM}, Th. 1.2).
\end{rmk}

\end{document}